\documentclass[twocolumn]{autart}    

\usepackage{graphicx}
\usepackage{hyperref}
\usepackage{amsmath,amssymb}
\usepackage{subcaption}
\usepackage{tikz}
\usepackage{todonotes}
\usetikzlibrary{shapes}
\usepackage[noend]{algpseudocode}
\usepackage{algorithm}[boxed]
\usepackage{cite}
\newtheorem{remark}{Remark}

\begin{document}

\begin{frontmatter}
\title{Gradient Free Cooperative Seeking of a Moving Source\thanksref{footnoteinfo}} 

\thanks[footnoteinfo]{This paper was not presented at any IFAC
meeting. Corresponding author E. Michael}

\author[UniMelb]{Elad Michael}\ead{eladm@student.unimelb.edu.au}, %
\author[UniMelb]{Chris Manzie}\ead{manziec@unimelb.edu.au},       %
\author[EPFL]{Tony A. Wood}\ead{tony.wood@epfl.ch},        
\author[Technion]{Daniel Zelazo}\ead{dzelazo@technion.ac.il},      
\author[ANU]{Iman Shames}\ead{iman.shames@anu.edu.au}            

\address[UniMelb]{The University of Melbourne, Parkville VIC 3010, Australia} 
\address[EPFL]{SYCAMORE Lab, Ecole Polytechnique Federale de Lausanne (EPFL),
Lausanne, Switzerland}
\address[Technion]{TECHNION, Haifa 32000, Israel}        
\address[ANU]{CIICADA LAB, ANU, Canberra ACT 0200, Australia}      


\begin{abstract}                          
In this paper, we consider the optimisation of a time varying scalar field by a network of agents with no gradient information. We propose a composite control law, blending extremum seeking with formation control in order to converge to the extrema faster by minimising the gradient estimation error. By formalising the relationship between the formation and the gradient estimation error, we provide a novel analysis to prove the convergence of the network to a bounded neighbourhood of the field's time varying extrema. We assume the time-varying field satisfies the Polyak-\L{}ojasiewicz inequality and the gradient is Lipschitz continuous at each iteration. Numerical studies and comparisons are provided to support the theoretical results.
\end{abstract}

\end{frontmatter}
\section{Introduction}
Localising the source of an unknown or uncertain scalar field has attracted significant attention in recent years. Extremum seeking can then be understood as driving the state of an agent or network of agents to the source, and maintaining a steady state in the neighbourhood of this optimal state in the unknown field. The widespread applications include internal combustion engine calibration~\cite{killingsworth2009hcci}, locating RF leakage\cite{al2012position}, optimising energy distribution~\cite{ye2016distributed}, and mobile sensor networks~\cite{stankovic2011distributed}. The main challenge in general is the approximation of the field, or a valid descent direction, with the additional challenge in the multi-agent case of coordinating the agents to improve the estimation. In this work, we consider discrete time extremum seeking, for the more classical continuous time extremum seeking problem see~\cite{tan2010extremum} and the references therein.

Extremum seeking with a single agent primarily uses ``dither'' or other motion patterns to estimate a descent direction. In~\cite{cochran2009nonholonomic, zhang2007extremum}, extremum seeking with a single agent is investigated relying only on the measurements of the scalar field, without usage of the agent's position. Both approaches use a sinusoidal dither signal to estimate the gradient of the unknown field. Using finite difference with previous measurements, tracking and estimation error bounds for the minima of a time-varying scalar field are derived in~\cite{shames2019online}, along with extensive numerical studies using a single agent. A hybrid controller is defined in~\cite{mayhew2007robust}, conducting a series of line minimisations to construct the descent direction, with stability and convergence results. In~\cite{durr2013lie,durr2017extremum}, the authors derive an extremum seeking controller using Lie bracket approximations of the field, however the approach is only applied to continuous time dynamics and in static fields.

Using a network of agents allows for a more robust estimate of the gradient, as the measurements are typically assumed to be simultaneous and thus unaffected by a time-varying source. In\cite{biyik2008gradient} a network is used with a single leader determining the estimated gradient, employing a zero mean dither signal, with the followers only keeping formation. The authors show that with a fast dither and slow formation keeping, the followers only track the gradient descent movement of the leader. A game theoretic approach is used in~\cite{stankovic2011distributed} to find equilibria of each agent's individual cost functions, using local extremum seeking controllers with sinusoidal dither. Using multiple ``leader'' agents and only inter-agent bearing measurements, the authors in~\cite{zhao2015bearing,zhao2015translational} stabilise a formation in arbitrary dimension with leaders following reference velocities or trajectories. In addition, using only bearing measurements allows for formation scaling and rotation. In multi-agent approaches, the set of measurements from each agent can be used to compute an estimated gradient, assuming a single sensor aboard each agent\cite{khong2014multi,ogren2004cooperative,skobeleva2018planar,vandermeulen2017discrete,vweza2015gradient}. All of these publications use some form of the simplex gradient\cite{regis2015calculus}, as do we in this paper. The controller design derived in~\cite{khong2014multi} uses a centralised extremum seeking controller, with access to all of the agents' measurements, which provides reference velocities to each of the agents. Convergence guarantees are provided for a variety of formation and extremum seeking methods satisfying their assumptions. A centralised controller is implemented in~\cite{ogren2004cooperative} to track the estimated gradient using least squares estimation and refined by Kalman filtering. The agents are tasked with formation keeping around a virtual leader, which climbs the gradient of the unknown field. However, the problem formulation only considers finite manoeuvres, and the formation may move extremely slowly. For networks of $3$ agents in $2$ dimensions a distributed control law with exponential convergence guarantees is investigated in\cite{skobeleva2018planar}. The agents in\cite{vandermeulen2017discrete} use a dynamic consensus algorithm to coordinate the gradient estimation, combined with a zero mean dither to construct a local gradient estimation. Finally, in a series of papers\cite{circular2015distributed,circular2013consensus,circular2011collaborative,circular2010source}, a group of unicycle agents performing distributed extremum seeking in circular formations is examined. The agents stabilise their formation and gradient estimate using a consensus algorithm, and performs well even with lossy communication and time-varying communication networks. The algorithm described in~\cite{circular2015distributed} is implemented in Section~\ref{sec:simul} to compare to the results derived in this paper.

Recently, extremum seeking for sources with dynamics has received some attention. In Section 1.2 of~\cite{ariyur2003real}, an extremum seeking algorithm using the internal model principle is derived, but requires extensive internal knowledge of the plant's dynamics. In~\cite{poveda2021fixed}, the authors derive similar tracking results to those provided here, albeit using continuous time dynamics and with the assumption of strong convexity using a non-smooth extremum seeking controller. Several recent works address time-varying extremum seeking in continuous time using a periodic dither algorithm~\cite{hazeleger2020extremum,grushkovskaya2017extremum,moshksar2015model}, however these works assume that the source/plant variation is significantly slower than the dither speed to allow for gradient estimation.

\paragraph*{Contributions}
This paper provides a novel analysis of multiagent extremum seeking focused on a time-varying source without using a centralised coordinator or dither motion, with discrete dynamics. This differs from the majority of the literature, which assumes a static or slowly drifting scalar field. In this work,
\begin{itemize}
    \item we allow the scalar field to be time-varying with no constraints on periodicity or time-scale separation;
    \item we incorporate formation control into extremum seeking using a novel condition on the formation potential, formalising the relationship between the gradient estimation and the formation;
    \item we show that the agents converge to a bounded neighbourhood of the time-varying extrema of the field;
    \item we present two elliptical error bounds on the gradient approximation of a function with a Lipschitz continuous gradient.
    \item we provide an open-source implementation of the approach to allow for further research and validation of our results.
\end{itemize}
Finally, at each iteration, we only assume that the time varying field is represented by a function which has Lipschitz continuous gradient (bounded second derivative), and satisfies the Polyak-\L{}ojasiewicz (PL) inequality. The PL inequality assumption is also weaker than many which are used to provide the linear convergence of gradient descent algorithms, such as convexity or quadratic growth\cite{plstuff}. The authors' previous investigation into this problem \cite{michael2020optimisation} included a more complicated control law than is presented here, with results restricted to $2$ dimensions. In this analysis, we simplify the control law, derive stronger convergence guarantees, and broaden the method to arbitrary dimension.

The paper is organised as follows. Section~\ref{sec:probForm} is devoted to basic assumptions on the time-varying field and agent dynamics. Section~\ref{sec:coopGradDesc} discusses the distributed control law and its performance for extremum seeking and formation keeping. Section~\ref{sec:gradEstimation} provides an example of cooperative gradient estimation, an improvement and generalisation of the results from~\cite{michael2020optimisation}. We provide numerical simulations in Section~\ref{sec:simul}, and conclude in Section~\ref{sec:conclude}.
\section{Problem Formulation}\label{sec:probForm}
Consider a network of $n$ agents where $x_k^{(i)} \in \mathbb{R}^d$ denotes the position of the $i$-th agent for $i\in \{1,...,n\}$ at iteration $k$. We use bold variables throughout the paper to describe the stacked vector for all agents, i.e. ${\bf x}_k$ to denote the vector of all agents' states stacked vertically. Let $\mathcal{G}=(\mathcal{V},\mathcal{E})$ be the underlying graph of the network with the vertex set $\mathcal{V}=\{1,...,n\}$ representing the agents and the edge set $\mathcal{E}\subseteq \mathcal{V}\times\mathcal{V}$ representing the communication topology. For each agent $i$, we define a set of neighbours $\mathcal{N}^{(i)}:= \{ j \mid (j,i)\in\mathcal{E} \}$ from which agent $i$ receives information at each iteration step.
\begin{assum}\label{ass:networkProp}
Assume that the agent communication graph $\mathcal{G}=(\mathcal{V},\mathcal{E})$ is connected and time invariant.
\end{assum}
The agents are modelled as single integrators:
\begin{align}
x_{k+1}^{(i)} = x_k^{(i)} + \alpha p_k^{(i)}, \label{eq:dyn}
\end{align}
where $\alpha$ is a constant.
\begin{remark}
We assume the single integrator dynamics~\eqref{eq:dyn} to focus on the time-varying field and highlight the extremum seeking algorithm used. However, the proposed approach may provide waypoints for a lower level controller, which navigates on a faster timescale until the waypoint is reached and the next measurement collected. The extension from single integrator dynamics to more complicated dynamics including velocity saturated models and nonholonomic models is discussed in~\cite{zhao2017defend}.
\end{remark}
At each iteration $k$, the time-varying field is represented by the function $f_k:\mathbb{R}^{d}\rightarrow\mathbb{R}$ with the non-empty minimiser set $\mathcal{X}^*_{f_k} := \textrm{argmin}_{x\in\mathbb{R}^{d}} f_k(x)$. The agents can only measure the function value at their location at each iteration, i.e. the value $f_k(x_k^{(i)})$. For any dimension $m\in\mathbb{Z}^+$ we define the distance between a point $x\in\mathbb{R}^{m}$ and a set $\mathcal{S}\subseteq\mathbb{R}^{m}$ as
\begin{align}
    d(x,\mathcal{S}) = \inf_{y\in\mathcal{S}} ||y-x||, \label{eq:pointSetDist}
\end{align}
where $||\cdot||$ is the Euclidean norm. Additionally, for a function $h:\mathcal{D}\rightarrow\mathbb{R}$ we will use the shorthand
\begin{align}
    h^* := \inf_{x\in\mathcal{D}} h(x), \label{eq:funcMinVal}
\end{align}
to represent the minimum value of that function.

\begin{assum}(Differentiability and Lipschitz Gradient):\label{ass:Lipschitz}
For all $k\geq0$, the functions $f_k:\mathbb{R}^d \rightarrow \mathbb{R}$ are at least once continuously differentiable. The gradients are $L_{f}-$Lipschitz continuous, i.e. there exists a positive scalar $L_{f}$ such that, for all $k \geq 0, x\in\mathbb{R}^d,\; y\in\mathbb{R}^d$, $$|| \nabla f_k(x) - \nabla f_k(y) || \leq L_{f}||x-y||, $$ or equivalently $$ f_k(y) \leq f_k(x) + \nabla f_k(x)^T(y-x) + \frac{L_{f}}{2}||y-x||^2.$$
\end{assum}

\begin{assum}(Polyak-\L{}ojasiewicz (PL)Condition):\label{ass:polyak}
For all $k\geq0$, there exists a positive scalar $\mu_{f}$ such that $\frac{1}{2}||\nabla f_k(x)||^2\geq \mu_{f}(f_k(x) - f^*_k)$.
\end{assum}
The assumption that a function has an $L-$Lipschitz continuous gradient is equivalent to assuming the second derivative has bounded norm, if it is twice differentiable. The PL condition requires that the gradient grows faster than a quadratic as we move away from the optimal function value. The PL condition does not require the minima to be unique, although it does guarantee that every stationary point is a global minimum~\cite{plstuff}. In addition to Assumptions~\ref{ass:Lipschitz}-\ref{ass:polyak} on each $f_k$, we quantify the ``speed'' with which the field may vary next.
\begin{assum}(Bounded Drift in Time):\label{ass:drift}
There exist positive scalars $\eta_0$ and $\eta^*$ such that $|f_{k+1}(x)-f_k(x)|\leq \eta_0$ for all $x\in\mathbb{R}^d$ and $|f^*_k-f^*_{k+1}|\leq \eta^*$.
\end{assum}
The problem of interest is given below.
\begin{prob}\label{prob:onlyProb}
For a network of $n$ agents with dynamics~\eqref{eq:dyn} and communication topology satisfying Assumption~\ref{ass:networkProp}, let $\{f_k\}$ be a sequence of functions with a corresponding sequence of minimiser sets $\{\mathcal{X}^*_{f_k}\}$ satisfying Assumptions~\ref{ass:Lipschitz}-\ref{ass:drift}. Given the measurements $\mathcal{Y}_k^{(i)}=\{f_k(x_k^{(j)}) \mid j\in\mathcal{N}^{(i)}\cup\{i\}\}$, find $\alpha,p_k^{(i)}$ and a constant $M$ for all agents $i\in\mathcal{V}$ and for all $k\geq0$ such that $\lim\limits_{k\rightarrow\infty}d(x_k^{(i)},\mathcal{X}^*_{f_k})\leq M$.
\end{prob}

In Section~\ref{sec:coopGradDesc}, we will incorporate formation control into the extremum seeking algorithm. To this end, we use a formation potential function $\phi({\bf x}_k):\mathbb{R}^{nd}\rightarrow\mathbb{R}^+$ which takes the full state vector of all agents and returns a scalar which is minimised when the agents are in formation. Let the minimum be denoted by $\phi^* := \min\limits_{{\bf x}\in\mathbb{R}^{nd}}\phi({\bf x})$.
\begin{defn}\label{def:formFuncs}
We define $\phi:\mathbb{R}^{nd} \rightarrow \mathbb{R}^+$ to be the \emph{formation potential function} for the network, with minimisers $\mathcal{X}^*_{\phi}$, and assume the following properties. The function $\phi({\bf x}_k)$
\begin{enumerate}
    \item is continuously differentiable on $\mathbb{R}^{nd}$ with gradient which is Lipschitz continuous with constant $L_{\phi}$;
    \item satisfies the PL inequality (Assumption~\ref{ass:polyak}), with constant $\mu_{\phi} \geq \mu_{f}$;
    \item has gradient component $\nabla_{x^{(i)}_k} \phi({\bf x}_k)$ which is computable using only the state of agent $i$ and neighbours $j\in\mathcal{N}_{i}$;
\end{enumerate}
\end{defn}
In the definition of the formation potential functions, the first two conditions ensure that $\phi({\bf x}_k)$ shares the minimal properties that make $f_k$ amenable to analysis. The third property ensures that the local information each agent has is sufficient for computation of the descent direction. Navigation potential functions such as in~\cite{do2006formation,tanner2005formation,olfati2002distributed,de2006decentralized,dimarogonas2006feedback} satisfy these assumptions, and we give a further simple example of a formation potential function which satisfies these assumptions in Section~\ref{sec:simul}.

\section{Cooperative Gradient Descent}\label{sec:coopGradDesc}

In this section we discuss our primary result, showing that a network of agents cooperating can reach a bounded neighbourhood of the minimiser set. In this section, for simplicity, we assume each agent uses an $\varepsilon-$\emph{gradient oracle} at each iteration to construct a step direction.

\begin{defn}\label{def:epsOracle}
\emph{($\varepsilon$-gradient oracle):} Given the function $f_k:\mathbb{R}^{d}\rightarrow\mathbb{R}$ and the state of the agents in the network ${\bf x}_k\in\mathbb{R}^{d}$, the oracle returns $O(f_k,{\bf x}_k,\mathcal{N}^{(i)}) = \nabla f_k(x^{(i)}_k) + \varepsilon_k$.
\end{defn}

In order to motivate the incorporation of formation control, consider the case where $p^{(i)}_k = -O(f_k,{\bf x}_k,\mathcal{N}^{(i)})$:
\begin{align}
\begin{split}
x^{(i)}_{k+1} :=& x^{(i)}_{k} + \alpha (-O(f_k,{\bf x}_k,\mathcal{N}^{(i)})) \\
=& x^{(i)}_{k} - \alpha (\nabla f_k(x^{(i)}_k) + \varepsilon_k)),
\end{split}\label{eq:naiveDyn}
\end{align}
and provide the following lemma on the convergence properties of the system.

\begin{lem}\label{lem:noisyGradDesc}
For a sequence of functions $\{f_k\}$ with minimiser sets $\{\mathcal{X}^*_{f_k}\}$ satisfying Assumptions~\ref{ass:Lipschitz}-\ref{ass:drift}, the system with dynamics~\eqref{eq:naiveDyn} satisfies
\begin{align}
\begin{split}
\frac{1}{2}d(x^{(i)}_k,\mathcal{X}^*_{f_k})^2 &\leq \beta(d(x^{(i)}_0,\mathcal{X}^*_{f_0})^2,k) \\
&\hspace{-1cm} + \frac{\alpha}{2\mu_{f}}\sum_{t=0}^{k}(1 - \alpha\mu_{f})^{k-t}||\varepsilon_t||^2  + \frac{\eta_{0}+\eta^*}{\alpha\mu_{f}^2}, \label{eq:noisyGradNeighbourhood}
\end{split}
\end{align}
for $\beta\in\mathcal{KL}$, $\alpha\in(0,\frac{1}{L_{f}}]$ with $L_{f},\mu_{f}$ from Assumptions~\ref{ass:Lipschitz}-\ref{ass:polyak}, and $d(x^{(i)}_k,\mathcal{X}^*_{f_k})$ defined in~\eqref{eq:pointSetDist}.
\end{lem}

\begin{pf}
See Appendix~\ref{app:noisyGradDesc}.
\end{pf}
\begin{rem}
  Lemma~\ref{lem:noisyGradDesc} seems to imply that if $\alpha$ is chosen to be $\frac{1}{\mu_f}$, the impact of the gradient error from steps before $k$ is zero. To understand why, note that the Lipschitz constant $L_{f}$ and PL constant $\mu_{f}$ satisfy the following
  \begin{align}
  \frac{\mu_{f}}{2}d(x^{(i)}_k,\mathcal{X}^*_{f_k})^2 \leq f_k(x) - f^*_k \leq \frac{L_{f}}{2}d(x^{(i)}_k,\mathcal{X}^*_{f_k})^2, \label{eq:quadraticSqueeze}
  \end{align}
  see~\cite{plstuff} for in depth discussion regarding the PL inequality. Requiring that $\alpha \leq \frac{1}{L_{f}}$ implies $\alpha \leq \frac{1}{\mu_{f}}$. Thus, if $\alpha \approx \frac{1}{\mu_f}$, then we must have that $\mu_f\approx L_{f}$ and $f_k$ is approximately a scaled norm as a consequence of~\eqref{eq:quadraticSqueeze}. For the scaled norm function, the gradient dynamics~\eqref{eq:naiveDyn} would take the agent directly to the minimiser, except for the error term from the most recent gradient estimate in~\eqref{eq:noisyGradNeighbourhood} and the drift error term $\frac{\eta_0+\eta^*}{\mu_f}$.
\end{rem}
From Lemma~\ref{lem:noisyGradDesc}, the system with dynamics~\eqref{eq:naiveDyn} converges to a neighbourhood dependent on the magnitude of the gradient error terms $||\varepsilon_k||^2$ and a constant term due to drift. This result is similar to the $(\beta,\gamma)$-tracking property defined in~\cite{poveda2021fixed} with the power series in terms of $||\varepsilon_t||^2$ as the $\gamma$ function therein. As noted in their paper, the result from Lemma~\ref{lem:noisyGradDesc} resembles a semi-global practical ISS bound with respect to $||\varepsilon_k||^2$ as the input.
However, we can improve upon this tracking result, as the magnitude of the gradient error $||\varepsilon_t||^2$ is not bounded. An idea behind this work is that in using function samples to estimate the gradient, the error in estimation is generally a function of the geometry of the samples taken. By incorporating formation control into the dynamics, we are able to bound the error terms $||\varepsilon_k||^2$. We show a specific example of this in Section~\ref{sec:gradEstimation}, but make minimal assumptions in this section on the specifics of how to construct a gradient estimate from sample points.

To characterise the entire network's behaviour, we define the time-varying function $F_k:\mathbb{R}^{nd}\rightarrow\mathbb{R}$ 
\begin{align*}
F_k({\bf x}_k) := \sum_{i\in\mathcal{V}} f_k(x^{(i)}_k),
\end{align*}
and note that $F_k$ satisfies Assumptions~\ref{ass:Lipschitz}-\ref{ass:polyak} with the same constants $L_{f},\mu_{f}$. The time-varying minimiser set of $F_k({\bf x}_k)$ is
\begin{align*}
\mathcal{X}^*_{F_k} = \overbrace{\mathcal{X}^*_{f_k}\times\mathcal{X}^*_{f_k}...\times\mathcal{X}^*_{f_k}}^{n}.
\end{align*}

To incorporate formation control into the extremum seeking analysis, we make the following assumption about the selection of $\phi({\bf x}_k)$.
\begin{assum}\label{ass:formationAsError}
The formation potential function $\phi:\mathbb{R}^{nd} \rightarrow \mathbb{R}^+$ as in Definition~\ref{def:formFuncs} satisfies
\begin{align}
\phi({\bf x}_k) \geq \frac{c}{2}\sum_{i\in\mathcal{V}}||\varepsilon^{(i)}_k||^2,
\end{align}
where $\varepsilon^{(i)}_k$ is defined in Definition~\ref{def:epsOracle}, and $c\in\mathbb{R}^{+}$ chosen such that $c > \frac{1}{\mu_{f}}$.
\end{assum}

Assumption~\ref{ass:formationAsError} formalises the relationship between the gradient estimation error and the formation. In Section~\ref{sec:simul} we provide the example $\phi({\bf x}_k) = \phi^* + L_{f} \sum_{i\in\mathcal{V}} ||x^{(i)}-x^{(j)}-\hat{x}^{(ij)}||^2$, where the terms $\hat{x}^{(ij)}$ define the optimal formation and the constant $\phi^*$ ensures Assumption~\ref{ass:formationAsError} is satisfied when the agents are in perfect formation. The constant offset does not change the dynamics, it allows $\phi({\bf x}_k)$ to bound the gradient error in the convergence analysis, see the proof of Theorem~\ref{thm:compositeConvergence}. This formation potential function satisfies the assumptions in Definition~\ref{def:formFuncs} and Assumption~\ref{ass:formationAsError}, however it requires an \textit{apriori} selection of each agent's neighbours.

With the formation potential function defined, we define the ``composite'' function $\hat{f}_k:\mathbb{R}^{nd} \rightarrow \mathbb{R}$ as
\begin{align}
\hat{f}_k({\bf x}_k) := F_k({\bf x}_k) + \phi({\bf x}_k),\label{eq:compositeFunc}
\end{align}
with corresponding minimisers in the set $\mathcal{X}^*_{\hat{f}_k}$, and the new system dynamics
\begin{align}
x^{(i)}_{k+1} := x^{(i)}_k - \alpha(\nabla_{x^{(i)}_k} \hat{f}_k + \varepsilon_k). \label{eq:compositeDyn}
\end{align}
Each agent can compute the gradient $\nabla_{x^{(i)}_k} \phi({\bf x}_k)$ with only local information, so the gradient of the composite function, being the sum of $f_k$ and $\phi$, can be estimated by using the same $\varepsilon-$gradient oracle for $f_k$. Both $F_k$ and $\phi$ satisfy Assumption~\ref{ass:Lipschitz} with constants $L_{f},L_{\phi}$ respectively, and Assumption~\ref{ass:polyak} with constants $\mu_{f},\mu_{\phi}$. Therefore, the composite function satisfies both Assumptions~\ref{ass:Lipschitz}-\ref{ass:polyak} with constants $L_{\hat{f}} := L_{f}+L_{\phi}$ and $\mu_{\hat{f}} \geq \min(\mu_{f},\mu_{\phi}) = \mu_{f}$.

\begin{lem}\label{lem:minimiserRelation}
For the composite function $\hat{f}_k$, as defined in~\eqref{eq:compositeFunc}, we have
\begin{align*}
\hat{f}^*_k := \min_{{\bf x}\in\mathbb{R}^{nd}}\hat{f}_k({\bf x}) \leq \phi^* + \frac{\min(L_{f},L_{\phi})}{2} d(\mathcal{X}^*_{F_k},\mathcal{X}^*_{\phi})^2,
\end{align*}
where we define the distance between the minimiser sets as
\begin{align*}
d(\mathcal{X}^*_{F_k},\mathcal{X}^*_{\phi}) := \min \{ ||x^*_{\phi} - x^*_{F_k}|| \mid  x^*_{\phi}\in\mathcal{X}^*_{\phi} \; , \;x^*_{F_k}\in \mathcal{X}^*_{F_k} \}.
\end{align*}
\end{lem}

\begin{pf}
See Appendix~\ref{app:minimiserRelation}.
\end{pf}
In the following theorem, we show that by incorporating a formation potential function, which bounds the gradient estimation error, the agents converge to a bounded neighbourhood of the time varying minimiser set $\mathcal{X}^*_{\hat{f}_{k}}$. Furthermore, the system does not require leaders, a separate time-scale for the formation-keeping, or any centralised computation.

\begin{thm}\label{thm:compositeConvergence}
For a sequence of functions $\{\hat{f}_k\}$ as defined in~\eqref{eq:compositeFunc} with minimisers $\{\mathcal{X}^*_{\hat{f}_k}\}$, the system with dynamics~\eqref{eq:naiveDyn} satisfies
\begin{align*}
\begin{split}
\frac{1}{2}d(x^{(i)}_{k+1},\mathcal{X}^*_{\hat{f}_{k+1}})^2 &\leq \beta(d(x^{(i)}_{0},\mathcal{X}^*_{\hat{f}_0})^2,k) \\
&\hspace{-1cm} + \frac{\alpha}{c\mu}\sum_{t=0}^{k}(1 - \alpha\mu')^{k-t}\hat{f}^*_t + \frac{\eta_{0}+\eta^*}{\alpha\mu\mu'},
\end{split}
\end{align*}
for $\beta\in\mathcal{KL}$, $\alpha\in(0,\frac{1}{L_{\hat{f}}}]$, and $\mu' = \mu_{f}-\frac{1}{c}$. Therefore, we have
\begin{align}
\lim_{k\rightarrow\infty} \frac{1}{2}d(x^{(i)}_{k+1},\mathcal{X}^*_{\hat{f}_{k+1}})^2 \leq  \frac{\underset{k\rightarrow\infty}{\lim} \sup \hat{f}^*_k}{\mu'} + \frac{\eta_{0}+\eta^*}{\alpha\mu\mu'}. \label{eq:limBound}
\end{align}
\end{thm}

\begin{pf}
See Appendix~\ref{app:compositeConvergence}.
\end{pf}
\section{Gradient Estimation and Error}\label{sec:gradEstimation}
In Section~\ref{sec:coopGradDesc}, we assume that each agent has access to an estimate of $\nabla f(x_k^{(i)}) + \epsilon^{(i)}$. In this section, we provide a method by which agent $i$ can estimate $\nabla f(x_k^{(i)})$ as well as compute an error bound for the estimate. The error bound and gradient estimation method apply to \emph{any} function which satisfies Assumption~\ref{ass:Lipschitz}. This method is a significant improvement of our previous work\cite{michael2020optimisation} and generalises to any dimension with any number of neighbours. Furthermore, we emphasise that the results from Section~\ref{sec:coopGradDesc} are independent of this section. The results presented here are an example of one possible method of gradient estimation and estimation error bounding. We make the following assumption on the neighbour set.
\begin{assum}\label{ass:fullRankNeighbours}
For each agent $i\in\mathcal{V}$ with state $x^{(i)}_k\in\mathbb{R}^d$, the neighbour set cardinality satisfies $|\mathcal{N}^{(i)}| \geq d$. Further, the vectors $\{x^{(l)}_k - x^{(i)}_k\}_{l\in\mathcal{N}^{(i)}}$ span $\mathbb{R}^d$.
\end{assum}
The requirement that the agents do not arrange on a low dimensional subspace is one of the primary motivators for incorporating formation control, as well as preventing collisions in applications with physical robots. Similar requirements for the arranging of agents, and controllers to achieve non-collinearity, are discussed in~\cite{liu2021orthogonal,shames2013doppler,ogren2004cooperative,bishop2009bearing}.
\begin{remark}
In the absence of Assumption~\ref{ass:fullRankNeighbours}, it is still possible to compute an approximate gradient using a variety of methods, such as in~\eqref{eq:centreDef}. However, it is not possible to bound the error of the gradient estimate.
\end{remark}
We define three useful variables before proceeding:
\begin{align}
\begin{split}
s^{(ij)}_k &:= \frac{f_k(x^{(j)}_k) - f_k(x^{(i)}_k)}{||x^{(j)}_k - x^{(i)}_k||}, \\
v^{(ij)}_{k} &:= \frac{x^{(j)}_k - x^{(i)}_k}{||x^{(j)}_k - x^{(i)}_k||}, \\
a^{(ij)}_k &:= \frac{L_{f}}{2}||x^{(j)}_k - x^{(i)}_k||.
\end{split}
\label{eq:usefulConstants}
\end{align}
We use ${\bf s}^{(i)}_k,{\bf a}^{(i)}_k$ to denote the vertically stacked vectors of $s^{(ij)}_k,a^{(ij)}_k$ for all neighbours $j\in\mathcal{N}^{(i)}$.

\begin{lem}\label{lem:gradPolyhedron}
For a function $f_k$ satisfying Assumption~\ref{ass:Lipschitz} and an agent $i$ with neighbour set $\mathcal{N}^{(i)}$ satisfying Assumption~\ref{ass:fullRankNeighbours}, there exists a \textbf{bounded} polyhedron
\begin{align}
\mathcal{P}^{(i)}_k := \{ x\in\mathbb{R}^d \mid \begin{bmatrix} A^{(i)}_k \\ -A^{(i)}_k \end{bmatrix} x \leq b^{(i)}_k  \} \label{eq:polyhedron}
\end{align}
such that $\nabla f(x^{(i)}_k) \in \mathcal{P}^{(i)}_k$, for $A^{(i)}_k\in\mathbb{R}^{|\mathcal{N}^{(i)}|\times d}$ and $b_k\in\mathbb{R}^{2|\mathcal{N}^{(i)}|\times d}$.
\end{lem}

\begin{pf}
See Appendix~\ref{app:gradPolyhedron}.
\end{pf}
From Lemma~\ref{lem:gradPolyhedron}, there exists a bounded space $\mathcal{P}^{(i)}_k$ within which the gradient $\nabla f(x^{(i)}_k)$ must exist. In~\cite{michael2020optimisation}, we restricted the error bound analysis to $2$ dimensions with $2$ neighbours. The same method is not computationally feasible in higher dimension, as it requires computation of the largest diagonal in the $d$-parallelotope, which has $2^{d-1}$ diagonals. Instead, define the following ellipse
\begin{align}
m^{(i)}_{k} &:= \sqrt{\sum_{j\in\mathcal{N}^{(i)}} (|s^{(ij)}_k-(g^{(i)}_k)^Tv^{(ij)}_k| + a^{(ij)}_k)^2} \label{eq:ellipseScaling}\\
\mathcal{E}^{(i)}_k &:= \left\{ x\in\mathbb{R}^{d} \mid \left\|\frac{A^{(i)}_k(x-g^{(i)}_k)}{m^{(i)}_k}\right\|^2 \leq 1 \right\}, \label{eq:ellipseDef}
\end{align}
with $g^{(i)}_k$ the centre of $\mathcal{E}^{(i)}_k$, $A^{(i)}_k$ the matrix defined in Lemma~\ref{lem:gradPolyhedron}, and $s^{(ij)}_k,v^{(ij)}_k,a^{(ij)}_k$ defined in~\eqref{eq:usefulConstants}. Define
\begin{align}
g^{(i)}_k := ((A^{(i)}_k)^TA^{(i)}_k)^{\dagger}(A^{(i)}_k)^T{\bf s}^{(i)}_k. \label{eq:centreDef}
\end{align}
Note that the centre of the ellipse will serve as the gradient estimate for agent $i$, and is equivalent to the simplex gradient\cite{regis2015calculus} of agent $i$ and its neighbours. We use the superscript ${\dagger}$ to denote the Moore Penrose pseudo-inverse, which is equivalent to the inverse when Assumption~\ref{ass:fullRankNeighbours} holds. In the following theorem we present an error bound that is valid in arbitrary dimension for any number of neighbours.
\begin{thm}\label{thm:boundingEllipse}
For a function $f_k$ satisfying Assumption~\ref{ass:Lipschitz} and an agent $i$ with neighbour set $\mathcal{N}^{(i)}$ satisfying Assumption~\ref{ass:fullRankNeighbours}, let $\mathcal{P}^{(i)}_k$ be the polytope defined in Lemma~\ref{lem:gradPolyhedron}.
Then $\mathcal{P}^{(i)}_k\subseteq \mathcal{E}^{(i)}_k$, for $\mathcal{E}^{(i)}_k$ the ellipse defined in~\eqref{eq:ellipseDef} with centre $g^{(i)}_k$ defined in~\eqref{eq:centreDef}.
Further, if $|\mathcal{N}^{(i)}| = d$, and we assume $B(r,c) = \{ x\in\mathbb{R}^{d} \mid ||x-c||_2 \leq r\}$ is the smallest bounding ball such that $\mathcal{P}^{(i)}_k \subseteq B(r,c)$, then
\begin{align}
\frac{||{\bf a}_k^{(i)}||}{\sigma_{\textrm{max}}(A^{(i)}_k)} \leq r \leq \frac{||{\bf a}_k^{(i)}||}{\sigma_{\textrm{min}}(A^{(i)}_k)} \label{eq:radiusBounds}
\end{align}
for $\sigma_{\textrm{max/min}}$ the largest/smallest singular values of $A^{(i)}_k$ and ${\bf a}_k^{(i)}$ the vector of $a_k^{(ij)}$ for all $j\in\mathcal{N}^{(i)}$.
\end{thm}

\begin{figure}[thpb]
 \centering
 \includegraphics[scale=0.45]{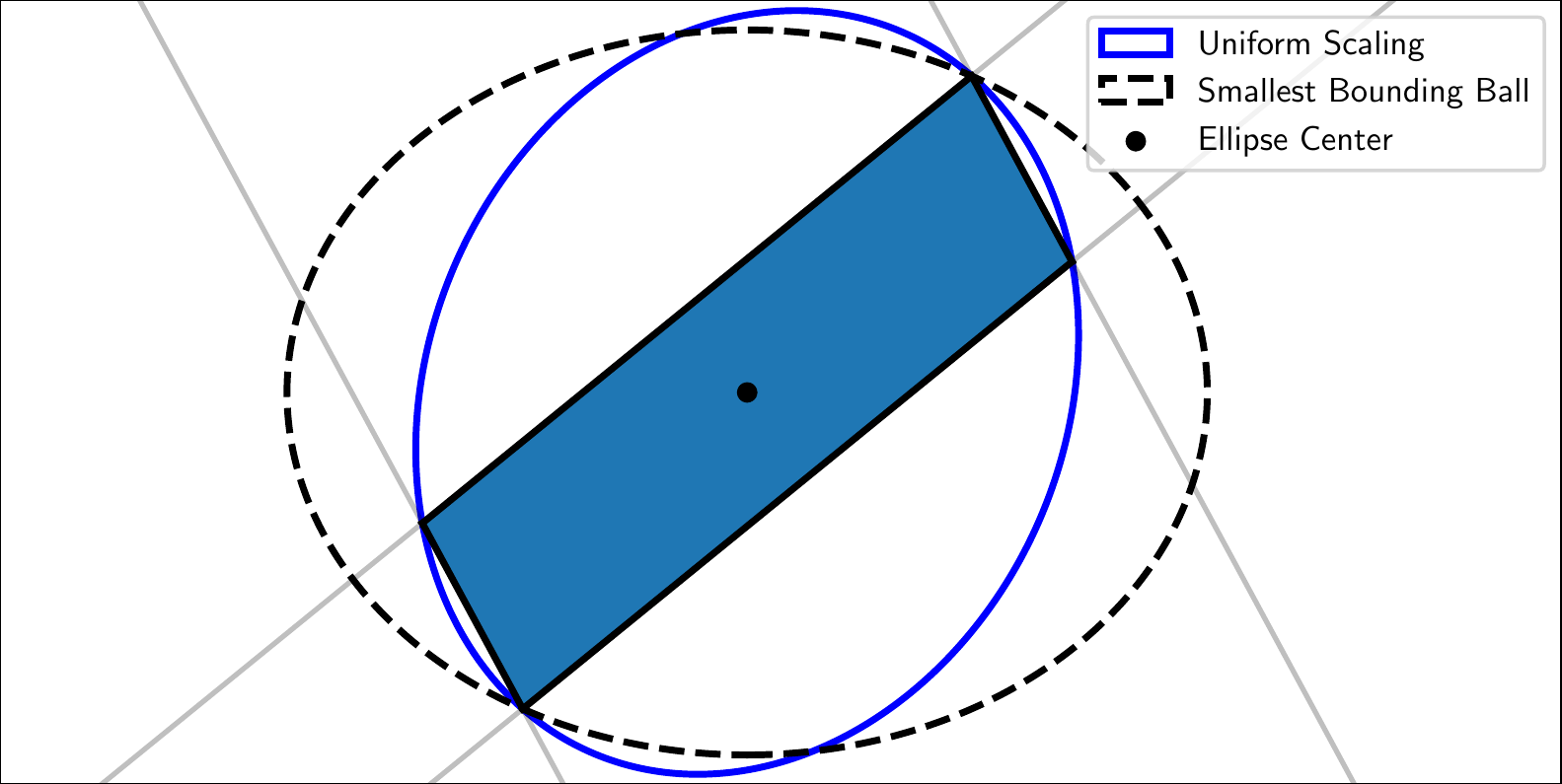}
 \caption{Ellipse bounding demonstration of Theorem~\ref{thm:boundingEllipse}. \label{fig:2neigh}}
\end{figure}

\begin{pf}
See Appendix~\ref{app:boundingEllipse}.
\end{pf}
The result in~\eqref{eq:radiusBounds} may be interpreted as ``the radius of the smallest bounding ball lies between the largest and smallest radii of $\mathcal{E}^{(i)}_k$.'' A simple example of the ellipse~\eqref{eq:ellipseDef} with $2$ neighbours labelled \emph{uniform scaling} (due to the uniform scaling of the shape matrix) compared the smallest bounding ball is shown in Figure~\ref{fig:2neigh}.

Given that finding the smallest bounding ball which contains a polytope is an NP hard problem, even for the relatively simple centrally symmetric parallelotopes\cite{bodlaender1990computational}, this approximation is sufficient for the primary goal of gradient estimation. Further, this approximation method gives the smallest $2$-norm bound on the error in the simplest case, with $d$ neighbours distributed in a lattice around agent $i$, as demonstrated in Corollary~\ref{cor:boundingBall}.
\begin{cor}\label{cor:boundingBall}
If agent $i$ has neighbour set with cardinality $|\mathcal{N}^{(i)}|=d$, and $(v^{(ij)}_k)^Tv^{(il)}_k=0$ for all $j,l\in\mathcal{N}^{(i)}$ with $j\neq l$, then $\mathcal{E}^{(i)}_k$ as defined in~\eqref{eq:ellipseDef} is the smallest bounding ball such that $\mathcal{P}^{(i)}_k\in\mathcal{E}^{(i)}_k$.
\end{cor}
\begin{pf}
If all neighbours are orthogonal, then $A^{(i)}_k$ as defined in Lemma~\ref{lem:gradPolyhedron} is an orthogonal matrix, i.e. $(A^{(i)}_k)^TA^{(i)}_k = I$. Therefore, $\mathcal{E}^{(i)}$ is a ball. Further, from Theorem~\ref{thm:boundingEllipse}, the smallest bounding ball radius lies between the largest and smallest radii of $\mathcal{E}^{(i)}_k$, which in this case are the same radius. Therefore, $\mathcal{E}^{(i)}$ is the smallest bounding ball containing $\mathcal{P}^{(i)}_k$.
\end{pf}
\vspace{-0.25cm}
For any number of neighbours satisfying Assumption~\ref{ass:fullRankNeighbours}, Theorem~\ref{thm:boundingEllipse} guarantees a gradient estimation error bound of the form
\begin{align}
||g^{(i)}_k - \nabla f_k(x^{(i)}_k)|| \leq \frac{m^{(i)}_k}{\sigma_{\min}(A_k^{(i)})}, \label{eq:gradErrorBound}
\end{align}
for $g^{(i)}_k$ the estimated gradient~\eqref{eq:centreDef} and $m^{(i)}_k$ as defined in~\eqref{eq:ellipseScaling}. Note that if Assumption~\ref{ass:fullRankNeighbours} does not hold, then $A_k^{(i)}$ is a low rank matrix, with a minimal singular value of $0$, and thus the bound~\eqref{eq:gradErrorBound} is undefined.

\subsection{Bounding Ellipse for large Neighbour Sets}
\begin{figure}[thpb]
 \centering
 \includegraphics[scale=0.5]{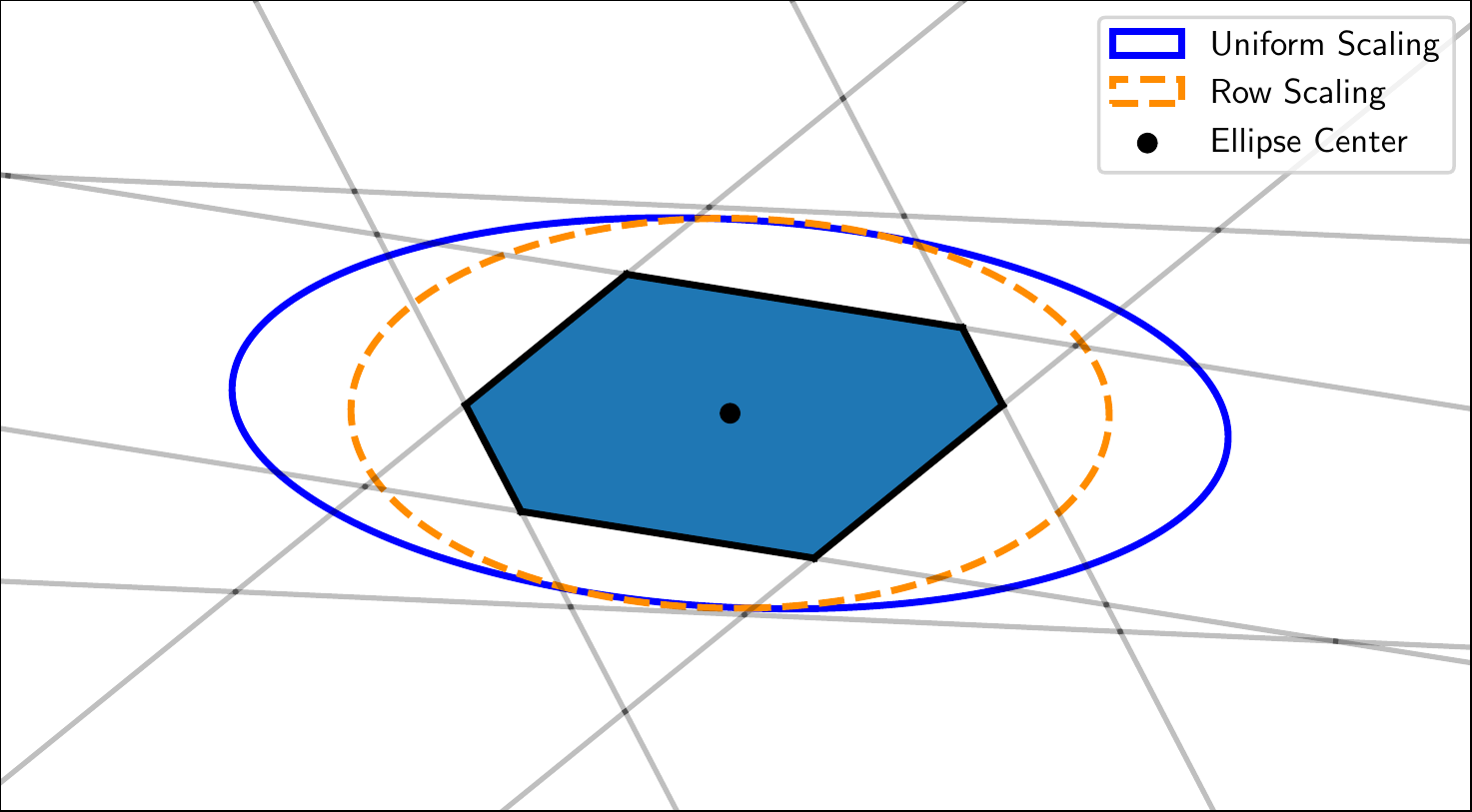}
 \caption{Comparing the bounds~\eqref{eq:ellipseDef} and~\eqref{eq:otherEllipseDef}. \label{fig:4neigh}}
\end{figure}

The ellipse from~\eqref{eq:ellipseDef} performs well for smaller sets of neighbours, but tends to be conservative when the neighbour set is larger than $d$. We provide an additional bounding ellipse here, which shares many of the useful properties of the ellipse defined in~\eqref{eq:ellipseDef}, but tends to be significantly less conservative in larger problems. The potentially large scaling factor in the denominator of \eqref{eq:ellipseDef} is distributed to each row, rather than applied uniformly, which mitigates some of the inflation from redundant neighbours. We define a matrix $B^{(i)}_{k}\in\mathbb{R}^{|\mathcal{N}^{(i)}|\times d}$ with the $j$-th row $B^{(i)}_{k}[j]$ defined as
\begin{align}
B^{(i)}_{k}[j] := \frac{(v^{(ij)}_{k})^T}{\sqrt{|\mathcal{N}^{(i)}|}(|s^{(ij)}_k-(g^{(i)}_k)^Tv^{(ij)}_{k}| + a^{(ij)}_k)} \label{eq:otherEllipseRow}
\end{align}
for $g^{(i)}_k\in\mathbb{R}^{d}$ the centre of the ellipse. The second ellipsoidal approximation of $\mathcal{P}^{(i)}_k$ can then be defined as
\begin{align}
\bar{\mathcal{E}}^{(i)}_k := \left\{ x\in\mathbb{R}^{d} \mid \left\|B^{(i)}_k(x-g^{(i)}_k)\right\|^2 \leq 1 \right\}. \label{eq:otherEllipseDef}
\end{align}
It can be verified that $\bar{\mathcal{E}}^{(i)}_k$ defined in~\eqref{eq:otherEllipseDef} also contains $\mathcal{P}^{(i)}_k$. However, the radius of the smallest bounding ball is not guaranteed to lie between the largest and smallest eigenvalues, and thus $\bar{\mathcal{E}}^{(i)}_k$ does not satisfy the claims of Corollary~\ref{cor:boundingBall}. For problems with larger sets of neighbours, the authors have empirically observed that $\bar{\mathcal{E}}^{(i)}_k$ seems to be a tighter approximation of $\mathcal{P}^{(i)}_k$. An example comparing the ``uniform scaling ellipse'' from~\eqref{eq:ellipseDef} to the ``row scaling ellipse'' from~\eqref{eq:otherEllipseDef} is included in Figure~\ref{fig:4neigh}.

\section{Simulations}\label{sec:simul}

In this section we provide numerical studies to illustrate the results from the previous sections, as well as comparison to another distributed extremum seeking algorithm. For the time varying scalar field, we use convex quadratic functions $f_k(x) = \frac{1}{2}(x-c(k))^TQ(x-c(k)) + \zeta^T(x-c(k)) + p$, for positive semi-definite $Q$. The values used in the following plots are
\begin{align*}
Q = \begin{bmatrix} 2.66 &-0.36 \\ -0.35 &1.74 \end{bmatrix} \; , \; &\zeta = [-1.28,4.66]^T \;,\; p = 6.26,\\
c(k) = 10\sin(\frac{\sqrt{2}k}{100}) &+ 10\sin(\frac{\sqrt{3}k}{100}) +\frac{k}{100},
\end{align*}
with $L_f,\mu_f$ the largest and smallest eigenvalues of $Q$ respectively. For the formation control function, we designate a set of neighbours for each agent $\mathcal{N}^{(i)}$ along with a corresponding set of ideal displacements $\hat{x}^{(ij)}$. The formation potential function is then
\begin{align}
\phi({\bf x}_k) = \phi^* + L_{f}\sum_{i\in\mathcal{V}}\sum_{j\in\mathcal{N}^{(i)}} ||x^{(i)} - x^{(j)} - \hat{x}^{(ij)}||^2_2.  \label{eq:formationPotential}
\end{align}
For other potential functions which satisfy the definitions used here, see~\cite{do2006formation,tanner2005formation,olfati2002distributed,de2006decentralized,dimarogonas2006feedback}.
In~\cite{michael2020optimisation} we derive the error bound on the gradient estimation in two dimensions, and show that the estimation error is proportional to the distance between the agents, with proportionality constant $L_{f}$, so the the Lipschitz constant $L_{f}$ and the minimum value $\phi^*$ in~\eqref{eq:formationPotential} ensure that $\phi({\bf x}_k) $ satisfies Assumption~\ref{ass:formationAsError}. The minimum value $\phi^*$ is chosen as an upper bound on the gradient approximation error when the agents are in perfect formation, derived from the gradient estimation error bounds in Theorem~\ref{thm:boundingEllipse}.

The simulated methods include the composite method derived in Section~\ref{sec:coopGradDesc} using two different formations, as well as the consensus for circular formations from~\cite{circular2015distributed} for comparison. For the composite method, as described in Section~\ref{sec:coopGradDesc}, we use the simplex gradient as the local gradient estimation method at each iteration (Algorithm~\ref{alg:compDyn}).

\begin{algorithm}
\caption{Distributed Composite Dynamics\label{alg:compDyn}}
\begin{algorithmic}
\For{$k=1,2,...$}
    \For{$i\in\{1,2,...,n\}$}
        \State $g^{(i)}_k = ((A^{(i)}_k)^TA^{(i)}_k)^{\dagger}(A^{(i)}_k)^T{\bf s}^{(i)}_k$  
    \EndFor
    \For{$i\in\{1,2,...,n\}$}
        \State $x^{(i)}_{k+1} = x^{(i)}_{k} - \frac{1}{L}(g^{(i)}_k+\nabla_{x^{(i)}_k} \phi({\bf x}_k))$
    \EndFor
\EndFor
\end{algorithmic}
\end{algorithm}

The circular formation controller is presented in Algorithm~\ref{alg:circDyn}, and is written exactly as in~\cite{circular2015distributed} accounting for the notation of this paper. The parameters used within Algorithm~\ref{alg:circDyn} are the same as used in the original paper~\cite{circular2015distributed}, in the example provided therein without noise. The radius of the formation $D=3$, the rotation velocity $\omega=1$, $\epsilon=0.5$ and $\alpha=1$. The consensus matrix used is of the same form as in~\cite{circular2015distributed}, for $6$ agents we have used
\begin{align*}
 P = \begin{bmatrix}
 0.5 & 0.25 & 0 & 0 & 0 & 0.25 \\
 0.25 & 0.5 & 0.25 & 0 & 0 & 0 \\
 0 & 0.25 & 0.5 & 0.25 & 0 & 0 \\
 0 & 0 & 0.25 & 0.5 & 0.25 & 0 \\
 0 & 0 & 0 & 0.25 & 0.5 & 0.25 \\
 0.25 & 0 & 0 & 0 & 0.25 & 0.5
 \end{bmatrix}.
\end{align*}
\begin{algorithm}
\caption{Circular Source Seeking\label{alg:circDyn}}
\begin{algorithmic}
\For{$i=1,...,n$}
    \State $h^{(i)}_0 = \tilde{g}^{(i)}_0=h^{(i)}_{-1}=c^{(i)}_{0}+f_{0}(x^{(i)}_0)(x^{(i)}_0-c^{(i)}_0)$
    \State $\phi^{(i)} = i\frac{2\pi}{n}$\;
\EndFor
\For{$k=1,2,...$}
    \For{$i=1,...,n$}
        \State $g^{(i)}_{k}=c^{(i)}_{k}+\frac{2}{D^2}f(x^{(i)}_{k})(x^{(i)}_{k}-c^{(i)}_{k})$
        \State $\tilde{g}^{(i)}_{k}=(1-{\bf \alpha})\tilde{g}^{(i)}_{k-1}+\alpha\tilde{g}^{(i)}_{k}$
        \State $\tilde{h}^{(i)}_{k}=h^{(i)}_{k-1}+\tilde{g}^{(i)}_{k-1}-\tilde{g}^{(i)}_{k-2}$
    \EndFor
    \State ${\bf h_{k}} = (P\otimes I_{2})(\bf \tilde{h}_{k})$
    \For{$i=1,...,n$}
        \State $c^{(i)}_{k} = (1-\varepsilon)c^{(i)}_{k-1}+\varepsilon h^{(i)}_{k}$
        \State $x^{(i)}_{k} = c^{(i)}_{k} + D R(\phi^{(i)}+\omega k)$
    \EndFor
\EndFor
\end{algorithmic}
\end{algorithm}

\begin{figure}[t!]
    \centering
    \begin{subfigure}[t]{0.2\textwidth}
        \begin{tikzpicture}
        \node[draw,minimum size=3cm,regular polygon,regular polygon sides=6] (a) {};

        \foreach \x in {1,2,...,6}
          \fill (a.corner \x) circle[radius=2pt];

        \end{tikzpicture}
        \caption{Hexagonal\label{fig:formation2d}}
    \end{subfigure}%
    ~
    \begin{subfigure}[t]{0.2\textwidth}
      \begin{tikzpicture}
      \node[draw,minimum size=3cm,regular polygon,regular polygon sides=4] (a) {};
      \node[draw,minimum size=3cm,regular polygon,regular polygon sides=4,right of= a,xshift=1.12cm] (b) {};

      \foreach \x in {1,2,...,4}
        \fill (a.corner \x) circle[radius=2pt];

      \foreach \x in {1,4}
        \fill (b.corner \x) circle[radius=2pt];

      \end{tikzpicture}
      \caption{Rectangular\label{fig:rectangleForm}}
    \end{subfigure}
    \caption{Neighbour topology for six agents in two dimensions.}
\end{figure}
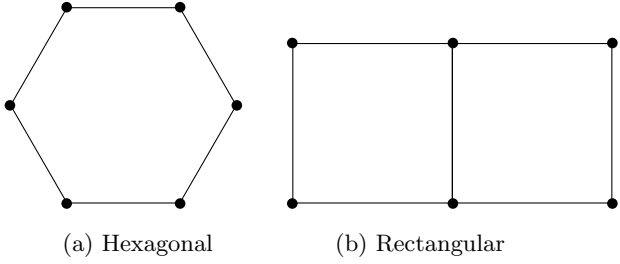

Choosing six agents forces the use of a regular hexagon for~\cite{circular2015distributed}. We therefore included the composite method using a regular hexagon formation for comparison. The neighbours are chosen to be the adjacent vertices as in Figure~\ref{fig:formation2d}.

While the circular motion controller in~\cite{circular2015distributed} requires this hexagonal arrangement for six agents, the framework proposed in this paper is flexible in the choice of formation by changing the ideal displacements $\hat{x}^{(ij)}$. To this end we also include a rectangular formation, illustrated in Figure~\ref{fig:rectangleForm}. As shown in~\cite{michael2020optimisation}, the gradient estimation error bound is a function of the orthogonality of the neighbours as well as the distance between them, so the rectangular formation will have lower gradient estimation error than the hexagonal formation with the same neighbour distances.
%
%
%

\begin{figure}[thpb]
 \centering
 \includegraphics[scale=0.45]{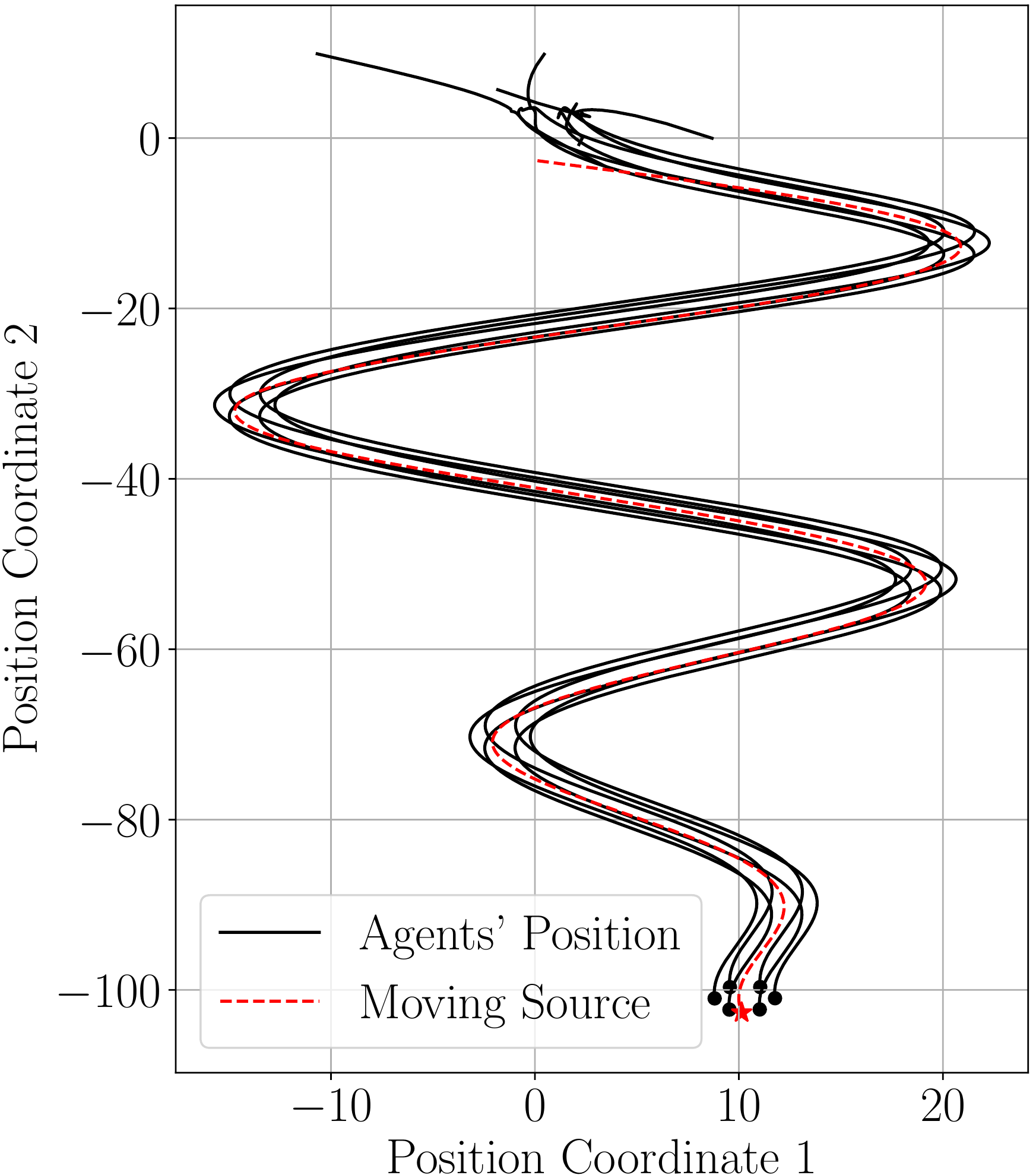}
 \caption{Agent Trajectories using the composite method from Section~\ref{sec:coopGradDesc}. \label{fig:trajComp}}
\end{figure}

\begin{figure*}[thpb]
 \centering
 \includegraphics[scale=0.33]{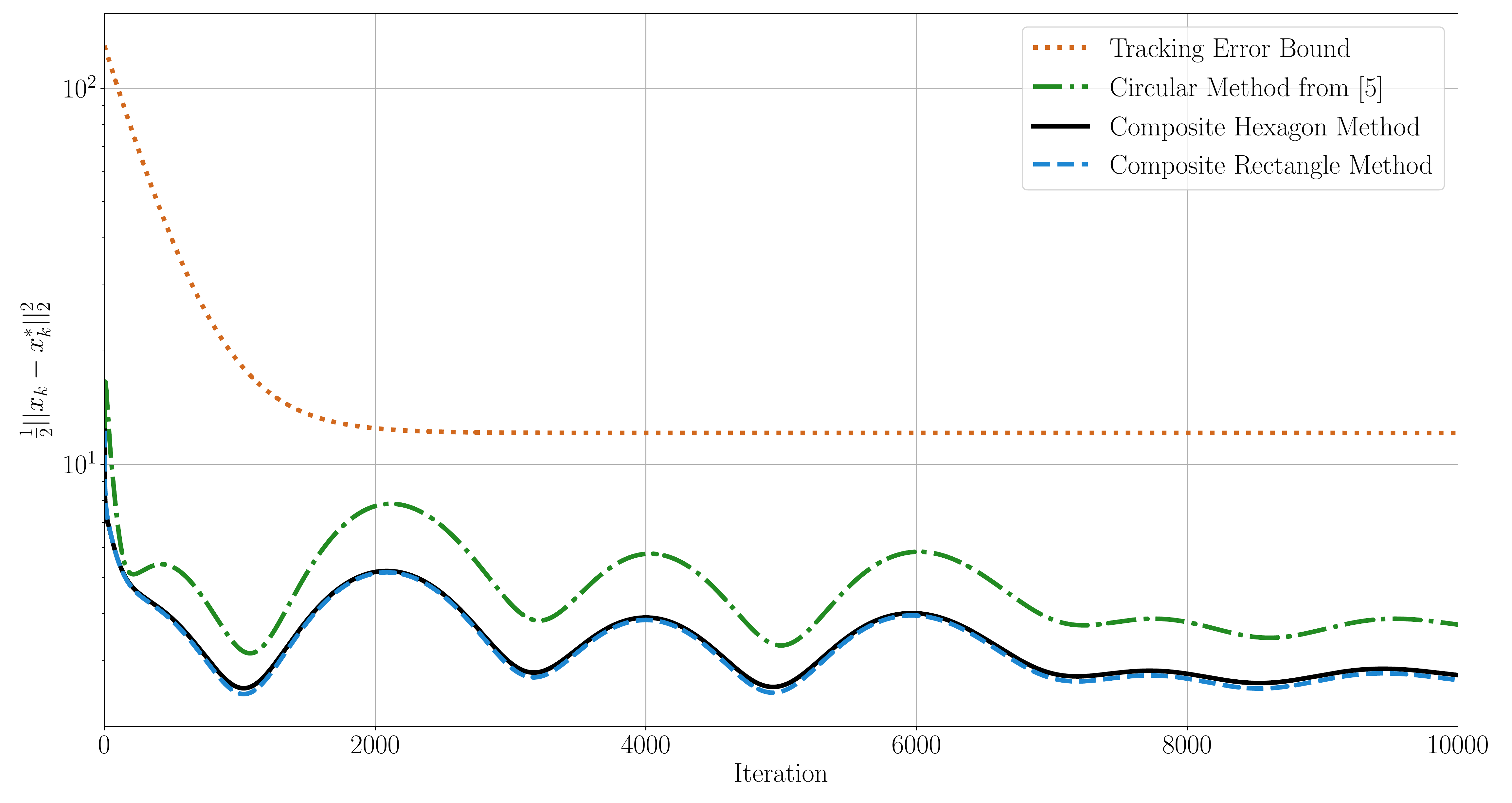}
 \caption{Comparison of formation distance from the signal source. \label{fig:minErrorComp}}
\end{figure*}

Figures~\ref{fig:trajComp} shows the resulting trajectories from the composite method. We exclude the trajectories from other methods, as they are visually identical. Instead, we include the comparison of the tracking error $\frac{1}{2}d(x_{k+1},\mathcal{X}^*_{\hat{f}_{k+1}})^2$ in Figure~\ref{fig:minErrorComp} for each method, including the theoretical bounds from Theorem~\ref{thm:compositeConvergence}.

We can see from Figure~\ref{fig:minErrorComp} that the theoretical minimiser error bound derived in Theorem~\ref{thm:compositeConvergence} holds in simulation. All methods exhibit similar performance, including the periodic increases in tracking error, i.e. the five ``bumps'' in Figure~\ref{fig:minErrorComp}. These coincide with the source accelerating around the curves of the path. The circular formation has higher tracking error, but the method in~\cite{circular2015distributed} is not explicitly designed to operate on time-varying scalar fields. The rectangular and hexagonal formations using the composite method track nearly identically, although the rectangular formation converges slightly closer to the optimal value set due to the lower gradient error.

\begin{figure}[thpb]
 \centering
 \includegraphics[scale=0.37]{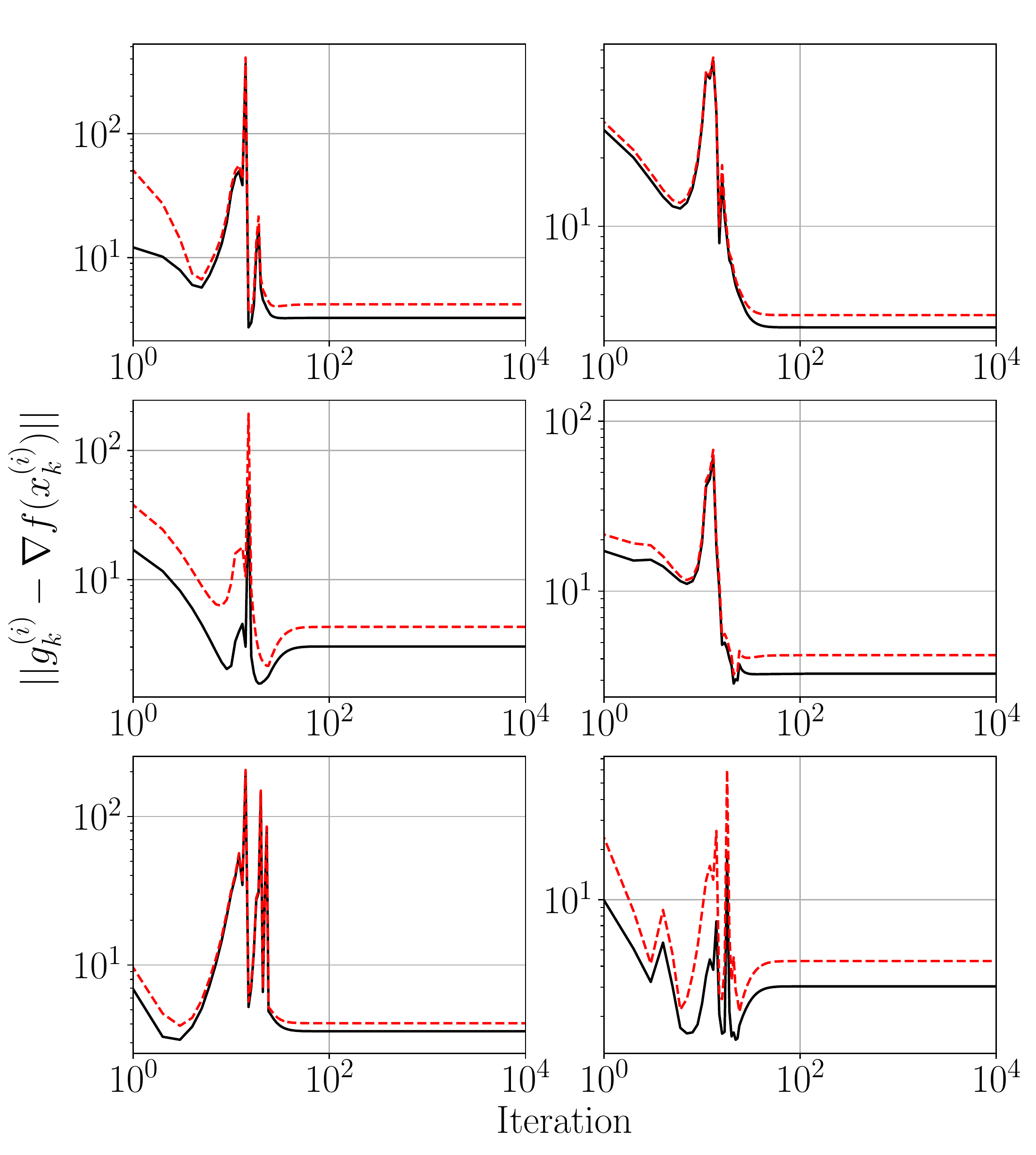}
 \caption{Gradient estimation error (black) and estimation error bound (red) for each agent from Figure~\ref{fig:trajComp} \label{fig:gradError}}
\end{figure}

In Figure~\ref{fig:gradError}, we show the error of the estimated gradient, as well as the error bound for each agent derived from the results of Theorem~\ref{thm:boundingEllipse}, defined in~\eqref{eq:gradErrorBound}.
Figure~\ref{fig:gradError} shows that once in formation, the agents gradient error (represented by the solid black line) becomes approximately constant. While gathering into formation early on however, there is a brief period of inaccurate gradient estimation. This further emphasises the importance of maintaining some formation to prevent collinearity. The gradient error bound (the dotted red line) tracks the variations of the gradient estimation error throughout the simulation, and is always within an order of magnitude or closer on this data set.

As the results from Section~\ref{sec:coopGradDesc} generalise to any dimension, we provide an example in three dimensions, as well as an implementation of the extremum seeking algorithm from Section~\ref{sec:coopGradDesc}, at the provided link.\footnote{\url{https://tinyurl.com/yc4fzpv2}}

\section{Conclusion}\label{sec:conclude}

In this paper we consider a formation of agents tracking the optimum of a time varying scalar field with no gradient information, in arbitrary dimension. At each iteration, the agents take measurements, communicate with their neighbours to estimate a descent direction, and converge to a neighbourhood of the optimum. We derive distributed control laws which drive the agents to a bounded neighbourhood of the optimiser set, without the delineation of leaders/followers or the use of communication intensive consensus protocols. The method is flexible to the choice of formation and gradient estimation method, and we provide examples using two formations and gradient estimation using the simplex gradient. By blending formation control with extremum seeking, the agents are able to minimise the gradient estimation error, improving the neighbourhood of convergence. We concluded with numerical studies showing that the proposed method is comparable with other extremum seeking methods, converging to a tighter neighbourhood while being more flexible in the choice of formation. Further research will focus on the relaxing of the assumptions on the formation potential functions, allowing for potential functions with non unique minima which do not satisfy the PL inequality, and incorporating time-varying neighbour sets.


\bibliographystyle{plain}        
\bibliography{references}           



\appendix

\section{Proof of Lemma~\ref{lem:noisyGradDesc}}\label{app:noisyGradDesc}

The agent identifying superscript $i$ is suppressed in this proof, as all calculations correspond to a single agent. By~\eqref{eq:naiveDyn} and Assumption~\ref{ass:Lipschitz} we have
\begin{align}
f_{k}(x_{k+1}) - f_k(x_k) &\leq \nabla f_k(x_k)^T (x_{k+1} - x_k) \nonumber\\
&\qquad + \frac{L_{f}}{2}||x_{k+1} - x_k||^2,\nonumber \\
&= - \alpha \nabla f_k(x_k)^T(\nabla f_k(x_k) + \varepsilon_k) \nonumber\\
&\qquad + \frac{\alpha^2L_{f}}{2}||\nabla f_k(x) + \varepsilon_k||^2. \nonumber\\
\intertext{Adding and subtracting $\frac{\alpha}{2}||\varepsilon_k||^2$ to complete the square, we have}
f_{k}(x_{k+1}) - f_k(x_k) &\leq \frac{\alpha}{2}||\varepsilon_k||^2 -\frac{\alpha}{2}||\nabla f_k(x_k)||^2 \nonumber\\
&\qquad + \frac{\alpha}{2}(\alpha L_{f} - 1)||\nabla f_k(x) + \varepsilon_k||^2. \nonumber \\
\intertext{Given $\alpha\in (0,\frac{1}{L_{f}}]$, we have $\alpha L_{f} - 1 \leq 0$,}
f_{k}(x_{k+1}) - f_k(x_k) &\leq -\frac{\alpha}{2}||\nabla f_k(x_k)||^2 + \frac{\alpha}{2}||\varepsilon_k||^2. \nonumber \\
\intertext{Using Polyak-\L{}ojasiewicz bounds (Assumption~\ref{ass:polyak}),}
f_{k}(x_{k+1}) - f_k(x_k) &\leq -\alpha\mu_{f}(f_k(x_k) - f^*_k) + \frac{\alpha}{2}||\varepsilon_k||^2. \nonumber \\
\intertext{Adding $f_{k+1}(x_{k+1})-f_{k}(x_{k+1})+f_k(x_k)-f_{k+1}^*$ to both sides, and using the scalar bounds from Assumption~\ref{ass:drift}}
f_{k+1}(x_{k+1})-f^*_{k+1} &\leq f_k(x_k) - f^*_{k+1} + f_{k+1}(x_{k+1}) \nonumber\\
&\qquad - f_{k}(x_{k+1}) - \alpha\mu_{f}(f_k(x_k) - f^*_k)\nonumber\\
&\qquad + \frac{\alpha}{2}||\varepsilon_k||^2, \nonumber \\
f_{k+1}(x_{k+1})-f^*_{k+1} &\leq f_k(x_k) - f^*_{k} - \alpha\mu_{f}(f_k(x_k) - f^*_k)\nonumber\\
&\qquad + \frac{\alpha}{2}||\varepsilon_k||^2 + \eta^* + \eta_0, \nonumber\\
\begin{split}
f_{k+1}(x_{k+1})-f^*_{k+1} &\leq (1 - \alpha\mu_{f})(f_k(x_k) - f^*_k)\\
&\qquad + \frac{\alpha}{2}||\varepsilon_k||^2 + \eta^* + \eta_0.
\end{split}\label{eq:recursive}
\end{align}
We can expand the recursive relationship~\eqref{eq:recursive} to the initial conditions
\begin{align}
\begin{split}
f_{k+1}(x_{k+1})-f^*_{k+1} &\leq (1 - \alpha\mu_{f})^{k}(f_0(x_0) - f^*_0)\\
& + \frac{\alpha}{2}\sum_{t=0}^{k}(1 - \alpha\mu_{f})^{k-t}||\varepsilon_t||^2\\
& + (\eta^* + \eta_0)\frac{1 - (1-\alpha\mu_{f})^k}{\mu_{f}\alpha}.
\end{split} \label{eq:subOptBndNaive}
\end{align}
Converting this sub-optimality bound into a bound on the convergence neighbourhood, we use the relationships~\eqref{eq:quadraticSqueeze} to obtain the final result
\begin{align*}
\frac{1}{2}d(x_{k+1},\mathcal{X}^*_{f_{k+1}})^2 &\leq \frac{(1 - \alpha\mu_{f})^{k}}{\mu_{f}}(\frac{L_{f}}{2}d(x_{0},\mathcal{X}^*_{f_{0}})^2 - \eta^* - \eta_0) \\
&+ \frac{\alpha}{2\mu_{f}}\sum_{t=0}^{k}(1 - \alpha\mu_{f})^{k-t}||\varepsilon_t||^2 + \frac{\eta^* + \eta_0}{\mu_{f}^2\alpha}.
\end{align*}

\section{Proof of Lemma~\ref{lem:minimiserRelation}}\label{app:minimiserRelation}

Let $x^*_{\phi}\in\mathcal{X}^*_{\phi}$ and $x^*_{F_k}\in\mathcal{X}^*_{F_k}$ be any of the points satisfying $||x^*_{\phi} - x^*_{F_k}|| = d(\mathcal{X}^*_{\phi},\mathcal{X}^*_{F_k})$. By the Lipschitz property of $F_k$ we have
\begin{align}
F_k(x^*_{\phi}) &\leq \frac{L_{f}}{2}||x^*_{F_k}-x^*_{\phi}||^2,\nonumber\\
& = \frac{L_{f}}{2} d(\mathcal{X}^*_{\phi},\mathcal{X}^*_{F_k})^2. \label{eq:compBound1}
\end{align}
By the Lipschitz property of $\phi$ we have
\begin{align}
\phi(x^*_{F_k}) &\leq \phi^* + \frac{L_{\phi}}{2}||x^*_{F_k}-x^*_{\phi}||^2,\nonumber\\
& = \phi^* + \frac{L_{\phi}}{2} d(\mathcal{X}^*_{\phi},\mathcal{X}^*_{F_k})^2. \label{eq:compBound2}
\end{align}
Using~\eqref{eq:compBound1}-\eqref{eq:compBound2}, we may bound the values of the composite function $\hat{f}_k$ at both $x^*_{\phi}$ and $x^*_{F_k}$. Therefore, given that the minimiser satisfies $\hat{f}^*_k \leq \hat{f}_k(x^*_{F_k})$ and $\hat{f}^*_k \leq \hat{f}_k(x^*_{\phi})$, we have
\begin{align*}
\hat{f}^*_k \leq \phi^* + \frac{\min(L_{f},L_{\phi})}{2} d(\mathcal{X}^*_{F_k},\mathcal{X}^*_{\phi})^2.
\end{align*}

\section{Proof of Theorem~\ref{thm:compositeConvergence}}\label{app:compositeConvergence}

Note that, just as in the proof of Lemma~\ref{lem:noisyGradDesc}, the agent identifying subscript is suppressed for readability as all calculations are with respect to one agent. As $\hat{f}_k$ shares all of the properties of $f_k$, we pick up from~\eqref{eq:recursive},
\begin{align*}
\hat{f}_{k+1}({\bf x}_{k+1})-\hat{f}^*_{k+1} &\leq (1 - \alpha\mu_{f})(\hat{f}_k({\bf x}_k) - \hat{f}^*_k)\\
&\qquad + \frac{\alpha}{2}||\varepsilon_k||^2 + \eta^* + \eta_0,
\intertext{Substituting the formation potential function $\frac{1}{c}\phi({\bf x}_k)$ for the error term $\frac{1}{2}||\varepsilon_k||^2$}
\hat{f}_{k+1}({\bf x}_{k+1})-\hat{f}^*_{k+1}  &\leq (1 - \alpha\mu_{f})(\hat{f}_k({\bf x}_k) - \hat{f}^*_k)\\
&\qquad + \frac{\alpha}{c}\phi({\bf x}_k) + \eta^* + \eta_0,
\intertext{Adding the strictly positive term $\frac{\alpha}{c} (f_k({\bf x}_k) - \hat{f}^*_k + \hat{f}^*_k)$ to the right side of the inequality}
\hat{f}_{k+1}({\bf x}_{k+1})-\hat{f}^*_{k+1} &\leq (1 - \alpha(\mu_{f}-\frac{1}{c}))(\hat{f}_k({\bf x}_k) - \hat{f}^*_k)\\
&\qquad + \frac{\alpha}{c}\hat{f}^*_k + \eta^* + \eta_0,
\end{align*}
Expanding the recursive relationship, with $\mu' := \mu_{f} - \frac{1}{c} \geq 0$, in terms of initial conditions yields
\begin{align}
\begin{split}
\hat{f}_{k+1}({\bf x}_{k+1})&-\hat{f}^*_{k+1} \leq \\
&(1 - \alpha\mu')^{k}(\hat{f}_0({\bf x_0}) - \hat{f}^*_0 - \eta^* - \eta_0)\\
&+ \frac{\alpha}{c}\sum_{t=0}^{k}(1 - \alpha\mu')^{k-t}\hat{f}^*_t + \frac{\eta^* + \eta_0}{\alpha\mu'}.
\end{split} \label{eq:subOptBndComplete}
\end{align}
Using~\eqref{eq:quadraticSqueeze}, as in the proof of Lemma~\ref{lem:noisyGradDesc}, we have
\begin{align*}
\frac{1}{2}d(x_{k+1},\mathcal{X}^*_{\hat{f}_{k+1}})^2 \leq &\frac{(1 - \alpha\mu')^{k}}{\mu_{f}}(\frac{L_{\hat{f}}}{2}d(x_{0},\mathcal{X}^*_{\hat{f}_{0}})^2 - \eta^* - \eta_0) \\
&+ \frac{\alpha}{c\mu_{f}}\sum_{t=0}^{k}(1 - \alpha\mu')^{k-t}\hat{f}^*_t + \frac{\eta^* + \eta_0}{\mu_{f}\mu'\alpha}.
\end{align*}

\section{Proof of Lemma~\ref{lem:gradPolyhedron}}\label{app:gradPolyhedron}

We begin by constructing the polyhedron $\mathcal{P}^{(i)}_k$, and showing that $\nabla f_k(x^{(i)}_k) \in \mathcal{P}^{(i)}_k$. We then show that, if Assumption~\ref{ass:fullRankNeighbours} holds, the polyhedron is bounded. None of the following analysis spans iterations, so we suppress the iteration counter $k$ for simplicity.

Consider agents $x^{(i)},x^{(j)},x^{(l)}\in\mathbb{R}^{d}$ with $j,l\in\mathcal{N}^{(i)}$. By the mean value theorem, we have that there exists a $t\in[0,1]$ such that
\begin{align}
\begin{split}
\hspace{-0.15cm}\nabla f((1-t)x^{(i)} + tx^{(j)})^T v^{(ij)} &= \frac{f(x^{(j)}) - f(x^{(i)})}{||x^{(j)}-x^{(i)}||}\\
& = s^{(ij)},
\end{split} \label{eq:meanValue}
\end{align}
for $v^{(ij)},s^{(ij)}$ defined in~\eqref{eq:usefulConstants}. On the right of~\eqref{eq:meanValue} we have the average directional derivative, which we will use to estimate the true directional derivative at $x^{(i)}$. Combining~\eqref{eq:meanValue} with Assumption~\ref{ass:Lipschitz} gives the worst case error of the directional derivative estimation,
\begin{align}
||\nabla f(x^{(i)})^T v^{(ij)} - s^{(ij)} || \leq \frac{L_{f}}{2}||x^{(ij)}|| = a^{(ij)}. \label{eq:directionalErrorBound}
\end{align}
We may rearrange~\eqref{eq:directionalErrorBound} into a pair of inequalities
\begin{align}
\begin{split}
(v^{(ij)})^T \nabla f(x^{(i)})  &\leq s^{(ij)} + a^{(ij)}\\
(-v^{(ij)})^T \nabla f(x^{(i)})  &\leq  a^{(ij)} - s^{(ij)}.
\end{split}\label{eq:boundingPlanes}
\end{align}
The two inequalities in~\eqref{eq:boundingPlanes} represent two hyperplanes within which the gradient is constrained. The two are oriented by the normal vector $v^{(ij)}$, separated by $2a^{(ij)}$, and centred on the plane $(v^{(ij)})^Tx = s^{(ij)}$. Define the matrix $A\in\mathbb{R}^{|\mathcal{N}^{(i)}|\times d}$, with each row equal to $v^{(ij)}$ for a neighbour $j\in\mathcal{N}^{(i)}$, and a vector $b\in\mathbb{R}^{2|\mathcal{N}^{(i)}|}$, with $s^{(ij)} + a^{(ij)}$ for each neighbour $j\in\mathcal{N}^{(i)}$ stacked above $a^{(ij)} - s^{(ij)}$ for each neighbour. Then the definition of the polyhedron $\mathcal{P}^{(i)}$ from Lemma~\ref{lem:gradPolyhedron} represents the set of $2|\mathcal{N}^{(i)}|$ inequalities from~\eqref{eq:boundingPlanes}, and we have $\nabla f(x^{(i)}) \in \mathcal{P}^{(i)}$.

To see that the polyhedron is bounded, let $\{e_1,e_2,...,e_n\}$ be the set of canonical basis vectors in $\mathbb{R}^d$. By Assumption~\ref{ass:fullRankNeighbours}, the vectors $\{v^{(ij)}\}_{j\in\mathcal{N}^{(i)}}$ span $\mathbb{R}^d$, and we may express each basis vector by a linear combination $e_l = \sum_{j\in\mathcal{N}^{(i)}} c^{(j)}_lv^{(ij)}$. We then have, for each point $x\in\mathcal{P}$,
\begin{align*}
e_l^Tx &= \sum_{j\in\mathcal{N}^{(i)}} (c^{(j)}_lv^{(ij)})^Tx \\
&\leq \sum_{j\in\mathcal{N}^{(i)}} c^{(j)}_l(s^{(ij)} + a^{(ij)})
\end{align*}
using the first inequality from~\eqref{eq:boundingPlanes}. In the negative $e_l$ direction we make use of the second inequality in~\eqref{eq:boundingPlanes},
\begin{align*}
-e_l^Tx &= \sum_{j\in\mathcal{N}^{(i)}} (c^{(j)}_l (-v^{(ij)}))^Tx \\
&\leq \sum_{j\in\mathcal{N}^{(i)}} c^{(j)}_l(a^{(ij)} - s^{(ij)}).
\end{align*}
We therefore have that, if Assumption~\ref{ass:fullRankNeighbours} holds, the polyhedron is bounded in $\mathbb{R}^{d}$.

\section{Proof of Theorem~\ref{thm:boundingEllipse}}\label{app:boundingEllipse}
Once again, we suppress the iteration identifying subscript $k$, as all the analysis takes place in a single iteration. Define a shifted coordinate system $y = x - g^{(i)}$, with the centre of the ellipse $g^{(i)}$ as the origin. The inequalities defining the interior of the polytope $\mathcal{P}^{(i)}$ from~\eqref{eq:boundingPlanes} then become
\begin{align}
\begin{split}
(v^{(ij)})^T y  &\leq s^{(ij)}-(g^{(i)})^Tv^{(ij)} + a^{(ij)}\\
(-v^{(ij)})^T y  &\leq  a^{(ij)} - (s^{(ij)}-(g^{(i)})^Tv^{(ij)}).
\end{split}\label{eq:boundingShiftedPlanes}
\end{align}

Let $y\in\mathcal{P}^{(i)}$ be any point within the polytope, i.e. it satisfies~\eqref{eq:boundingShiftedPlanes} for all $j\in\mathcal{N}^{(i)}$. Then one of the following two inequalities hold
\begin{align*}
(y^Tv^{(ij)})^2 &\leq (s^{(ij)}-(g^{(i)})^Tv^{(ij)} + a^{(ij)})^2 \\
(y^Tv^{(ij)})^2 &\leq (a^{(ij)}-(s^{(ij)}-(g^{(i)})^Tv^{(ij)}))^2,
\end{align*}
depending on the sign of $s^{(ij)}-(g^{(i)})^Tv^{(ij)}$. We may then use the single inequality
\begin{align}
(y^Tv^{(ij)})^2 \leq (|s^{(ij)}-(g^{(i)})^Tv^{(ij)}| + a^{(ij)})^2, \label{eq:intBound}
\end{align}
for any point $y\in\mathcal{P}^{(i)}$. Given the matrix $A^{(i)}$ as defined Lemma~\ref{lem:gradPolyhedron}, we have
\begin{align*}
y^T(A^{(i)})^TA^{(i)}y &= \sum_{j\in\mathcal{N}^{(i)}} y^Tv^{(ij)}(v^{(ij)})^Ty \\
&= \sum_{j\in\mathcal{N}^{(i)}} ((v^{(ij)})^Ty)^2.
\intertext{Assuming $y\in\mathcal{P}^{(i)}$ and applying~\eqref{eq:intBound}}
y^T(A^{(i)})^TA^{(i)}y &\leq \sum_{j\in\mathcal{N}^{(i)}} (|s^{(ij)}-(g^{(i)})^Tv^{(ij)}| + a^{(ij)})^2,\\
||A^{(i)}y||^2 &\leq (m^{(i)}_k)^2,
\end{align*}
for $m^{(i)}_k$ defined in~\eqref{eq:ellipseScaling}. Therefore, we have that each point in the polytope $\mathcal{P}^{(i)}$ is in the ellipse~\eqref{eq:ellipseDef}. Note that this works for any centre $g^{(i)}$, but the resulting ellipse will be differently sized depending on the choice of $g^{(i)}$.

We now assume that $|\mathcal{N}^{(i)}|=d$, and therefore $\mathcal{P}^{(i)}$ is a $d$-parallelotope, with parallel and congruent opposite faces. The centre of the parallelotope $c$ is the point $A^{(i)}c = {\bf s}^{(i)}$, for ${\bf s}^{(i)}$ the vector of $s^{(ij)}$ for all $j\in\mathcal{N}^{(i)}$. This is the point at which all diagonals intersect, and are bisected, and thus must be the centre of the smallest bounding ball $\mathcal{B}^{(i)}(r,c)$. We note this point is also returned by~\eqref{eq:centreDef}, therefore the ellipse $\mathcal{E}^{(i)}$ and the smallest bounding ball share the same centre. We may then assume, without loss of generality, that the parallelotope is centred at the origin. This further simplifies the definition of $\mathcal{E}^{(i)}$, as the term $\sum_{j\in\mathcal{N}^{(i)}} (|s^{(ij)}_k-(g^{(i)})^Tv^{(ij)}_k| + a^{(ij)}_k)^2 = ||{\bf a}^{(i)}||^2$ for ${\bf a}^{(i)}$ the vector of $a^{(ij)}$ for all $j\in\mathcal{N}^{(i)}$ as defined in~\eqref{eq:usefulConstants}. Let $V_k\in\mathbb{R}^{n\times n}$ be a diagonal matrix with $V_{ii} \in \{-1,1\}$. Then the vertices of $\mathcal{P}^{(i)}$ are the points
\begin{align*}
v_k = (A^{(i)})^{-1}V_k{\bf a}^{(i)} \; \forall k\in[1,2,3,...,2^{n}].
\end{align*}
The smallest bounding ball, by definition, includes all these vertices and therefore
\begin{align}
r &\geq \max_{V_k} ||(A^{(i)})^{-1}V_k{\bf a}^{(i)}||, \label{eq:maxOverVert}\\
&\geq \sigma_{\textrm{min}}((A^{(i)})^{-1})||V_k{\bf a}^{(i)}||, \label{eq:lazyBound}\\
&= \frac{||{\bf a}^{(i)}||}{\sigma_{\textrm{max}}(A^{(i)})}. \label{eq:radiusLowerBound}
\end{align}

Furthermore, we have that $\mathcal{B}(c,r)$ has smaller radius than the largest radius of the ellipse $\mathcal{E}^{(i)}$, or there would trivially exist a smaller bounding ball. The largest radius of $\mathcal{E}^{i}$ corresponds to the inverse of the smallest singular value of the shape matrix, i.e.
\begin{align*}
r \leq \frac{||{\bf a}^{(i)}||}{\sigma_{\textrm{min}}(A^{(i)})}.
\end{align*}
Combining these results, and reorganizing, gives the bounds from the statement of the theorem
\begin{align*}
\frac{||\bf a^{(i)}||}{\sigma_{\textrm{max}}(A^{(i)}_k)} \leq r \leq \frac{||\bf a^{(i)}||}{\sigma_{\textrm{min}}(A^{(i)}_k)}.
\end{align*}

\end{document}